\documentclass[11pt]{amsart}
\usepackage[T1]{fontenc}
\usepackage[latin5]{inputenc}
\usepackage{float}
\usepackage{amssymb}
\usepackage{amsfonts}
\usepackage[centertags]{amsmath}
\usepackage{amsthm}
\usepackage{newlfont}
 \makeatletter
\theoremstyle{plain}
\newtheorem{thm}{Theorem}[section]
\numberwithin{equation}{section}
\numberwithin{figure}{section}  
\theoremstyle{plain}

\theoremstyle{plain}
\newtheorem{cor}[thm]{Corollary} 
\theoremstyle{plain}
\newtheorem{defin}[thm]{Definition}
\theoremstyle{plain}
\newtheorem{lem}[thm]{Lemma} 
\theoremstyle{plain}

\begin{document}
\title{On the Double Sequence Space $\mathcal{H}_{\vartheta}$ as an Extension of Hahn Space $h$}
\author{Orhan Tu\v{g}}
\author{Eberhard Malkowsky}
\author{Vladimir Rako\v{c}evi\'{c}}
\author{Taja Yaying}
\address[Orhan Tu\v{g}(Corresponding Author)]{Department of Mathematics Education, Faculty of Education, Tishk
	International University, Erbil, Iraq}
\email{orhan.tug@tiu.edu.iq}
\address[Eberhard Malkowsky]{Department of Mathematics, State University of Novi Pazar, Novi Pazar, Serbia.}
\email{ema@pmf.ni.ac.rs}
\address[Vladimir Rako\v{c}evi\'{c}]{Department of Mathematics,
	Faculty of Sciences and Mathematics, University of Ni\v{s},
	Vi\v{s}egradska 33, 18000, Ni\v{s}, Serbia}
\email{vrakoc@sbb.rs}

\address[Taja Yaying]{Department of Mathematics, Dera Natung Government College, Itanagar 791113, India}
\email{tajayaying20@gmail.com}

\begin{abstract}
Double sequence spaces have become a significant area of research within functional analysis due to their applications in various branches of mathematics and mathematical physics. In this study, we investigate Hahn double sequence space denoted as $\mathcal{H}_{\vartheta}$, where  $\vartheta\in\{p,bp,r\}$, as an extension of the Hahn sequence space $h$. Our investigation begins with an analysis of several topological properties of $\mathcal{H}_{\vartheta}$, apart from a comprehensive analysis of the relationship between Hahn double sequences and some other classical double sequence spaces. The $\alpha-$dual, algebraic dual and $\beta(bp)-$dual, and $\gamma-$dual of the space $\mathcal{H}_{\vartheta}$ are detrmined. Furthermore, we define the determining set of $\mathcal{H}_{\vartheta}$ and we state the conditions concerning the characterization of four-dimensional (4D) matrix classes $(\mathcal{H}_{\vartheta},\lambda)$, where $\lambda=\{\mathcal{H}_{\vartheta},\mathcal{BV}, \mathcal{BV}_{\vartheta 0}, \mathcal{CS}_{\vartheta},\mathcal{CS}_{\vartheta 0},\mathcal{BS}\}$ and $(\mu,\mathcal{H}_{\vartheta})$, where $\mu=\{\mathcal{L}_u, \mathcal{C}_{\vartheta 0}, \mathcal{C}_{\vartheta},\mathcal{M}_{u}\}$. In conclusion, this research contributes non-standard investigation and various significant results into the space $\mathcal{H}_{\vartheta}$. The conducted results are deepen the understanding of the space $\mathcal{H}_{\vartheta}$ and open up new avenues for further research and applications in sequence space theory.
\end{abstract}\maketitle

\section{Introduction, Notations, and Preliminaries}
By the symbol $\Omega,$ we mean the set consisting of all complex-valued double sequences. That is
$$\Omega=\{x=(x_{jk})_{j,k\in \mathbb{N}_0}: x_{j,k}\in \mathbb{C}\},$$
where the symbol $\mathbb{C}$ is the set of all complex numbers, and $\mathbb{N}$ the set of all natural numbers. Any linear subspace of $\Omega$ is called a double sequence space. The classical double sequence spaces, namely $\mathcal{M}_u$, $\mathcal{C}_{p}$, $\mathcal{C}_{p0}$, $\mathcal{C}_{bp}$, $\mathcal{C}_{bp0}$, $\mathcal{C}_{r}$, $\mathcal{C}_{r0}$, and $\mathcal{L}_q$. The respective definitions are as follows: $\mathcal{M}_u$ represents the space of all bounded sequences, $\mathcal{C}_{p}$ contains sequences convergent in Pringsheim's sense, $\mathcal{C}_{p0}$ denotes the null space of $\mathcal{C}_{p}$, $\mathcal{C}_{bp}$ consists of sequences in $\mathcal{M}_u\cap\mathcal{C}_{p}$, $\mathcal{C}_{bp0}=\mathcal{M}_u\cap\mathcal{C}_{p0}$, $\mathcal{C}_{r}$ accommodates regular convergent sequences, $\mathcal{C}_{r0}$ denotes the null space of $\mathcal{C}_{r}$, and $\mathcal{L}_q$ encompasses $q-$absolutely summable sequences, where $1\leq q<\infty$. As a significant contribution by M\'{o}ricz \cite{Moricz}, it was shown that the spaces $\mathcal{M}_u$, $\mathcal{C}_{bp}$, $\mathcal{C}_{bp0}$, $\mathcal{C}_{r}$, and $\mathcal{C}_{r0}$ are Banach spaces when endowed with the supremum norm $\|.\|_{\infty}$.

A double sequence space $\mu$ is said to be a $DK-$space if it is a local convex space, and all the semi-norms defined on $\mu$ by $x=(x_{kl})\to|x_{kl}|$ are continuous. When such a $DK-$ space occupies a Fr\`{e}chet topology, then it is called an $FDK-$space. If we have a normed $FDK-$space, then it is referred to as a $BDK-$space. Three well known examples of  $BDK-$spaces are $\mathcal{M}_u$, $\mathcal{C}_{bp}$ and $\mathcal{C}_{r}$, where the norm is given by
\begin{eqnarray*}
\|x\|_{\infty}=\sup_{k,l\in\mathbb{N}}|x_{kl}|<\infty.
\end{eqnarray*}

The double series is represented by the notation $\sum_{k,l=0}^{\infty}x_{kl}$ and is considered convergent with respect to $\vartheta$ if its $mn^{th}-$partial sum $s_{mn}=\sum_{k,l=0}^{m,n}x_{kl}$ converges as $m,n\to\infty$, that is, there exists $s\in \mathbb{R}$ in which $\vartheta-\lim_{m,n\to\infty}s_{mn}=\sum_{k,l=0}^{\infty}x_{kl}=s$.

The collection of $\vartheta-$convergent series is represented as $\mathcal{CS}_{\vartheta}$,  where $\vartheta\in\{p,p0,bp,bp0,r,r0\}$. Moreover by $\mathcal{BS}$, $\mathcal{BV}$, and $\mathcal{BV}_{\vartheta 0}=\mathcal{BV}\cap\mathcal{C}_{\vartheta 0}$, we represent the set of all bounded series, set comprising double sequences with bounded variation, and the null space of $\mathcal{BV}$, respectively ( see \cite{OVM3,BAFB}).

Sever \cite{SY} defined the residual component of a double series as follow:
\begin{eqnarray*}
	r_{mn}=\sum_{k,l=m+1,0}^{\infty,n}x_{kl}+\sum_{k,l=0,n+1}^{m,\infty}x_{kl}+\sum_{k,l=m+1,n+1}^{\infty,\infty}x_{kl}
\end{eqnarray*}
One can observe here that if a series is convergent by means of $\vartheta$, then the residual component of the series must be a null sequence. 

Let $x^{[m,n]}$ denote the $(m,n)^{th}-$section  of a double sequence $x=(x_{kl}),$  expressed by

\begin{eqnarray*}
x^{[m,n]}=\sum_{k,l=1}^{m,n}x_{kl}e^{(kl)},
\end{eqnarray*}
where $e^{(kl)}=(e_{mn}^{(kl)})_{m,n=1}^{\infty}$ for every $k,l\in\mathbb{N}$ is defined  with $e_{mn}^{(mn)}=1$ and $e_{mn}^{(kl)}=0$ for $(k,l)\neq(m,n)$ which was defined by Zeltser \cite{ZM2}. Let $\mu$ be an arbitrary double sequence space and suppose that it contains $\varphi$, where $\varphi=span\{e^{(kl)}\}$. Zeltser \cite{ZM2} defined the space $\mu$ to be an $AK(\vartheta)-$space if each $x\in\mu$ has the $AK(\vartheta)-$property in the set $\mu$, i.e., for an arbitrary $x=(x_{kl})$ we have,
\begin{eqnarray*}
x=\vartheta-\lim_{m,n\to\infty}x^{[m,n]}=\vartheta-\lim_{m,n\to\infty}\sum_{k,l=1}^{m,n}x_{kl}e^{(kl)}=\sum_{k,l=1}^{\infty}x_{kl}e^{(kl)}.
\end{eqnarray*}
Moreover, in the same paper \cite{ZM2}, Zeltser introduced the double sequences $e^l, e_k$ and $e$ with the following definitions 
\begin{enumerate}
\item[(i)] $e^l=\sum_ke^{kl}$, where every element in l-th column is $1$, and the rest are $0.$
\item[(ii)] $e_k=\sum_le^{kl}$, where all terms of the k-th row are $1$, and the rest are $0$,
\item[(iii)] $e=\sum_{kl}e^{kl}$, where all terms are $1$,
\end{enumerate}
for all $k,l,m,n\in\mathbb{N}$.

Let $E$ be a double sequence space. Then we define the differentiated and integrated spaces of $E$, respectively, as follows:
\begin{eqnarray}
	&&\label{differentiated}dE:=\left\{x=(x_{kl})\in\Omega:\left\{\frac{1}{kl}x_{kl}\right\}_{k,l\in\mathbb{N}}\in E \right\},\\
	&&\label{integrated}\int E:=\left\{x=(x_{kl})\in\Omega:\left\{klx_{kl}\right\}_{k,l\in\mathbb{N}}\in E \right\}.
\end{eqnarray}

 Let $G=(g_{mnkl})$ be a 4D infinite matrix. Then $G:\lambda\to\mu$ defines a matrix transformation., where $m,n,k,l\in\mathbb{N}$. That is, $G$ maps $x\in \lambda$ into $Gx\in \mu$, where
\begin{eqnarray}\label{eq1.1}
(Gx)_{mn}=\vartheta-\sum_{k,l}g_{mnkl}x_{kl}
\end{eqnarray}
for all  $m,n\in\mathbb{N}$.

Now we define the sets  $\mathcal{P}_{\infty}$, $\mathcal{P}_{\vartheta}$ and $\mathcal{P}_{\vartheta 0}$ in order to consider them in the main results of the paper. 

\begin{eqnarray*}
&&\mathcal{P}_{\infty}:=\left\{x=(x_{kl})\in\Omega:\sup_{m,n\in\mathbb{N}}\frac{1}{mn}\left|\sum_{k,l=1}^{m,n}x_{kl}\right|<\infty\right\},\\
&&\mathcal{P}_{\vartheta}:=\left\{x=(x_{kl})\in\Omega:\vartheta-\lim_{m,n\to\infty}\frac{1}{mn}\sum_{k,l=1}^{m,n}x_{kl} \textit{ exists }\right\},\\
&&\mathcal{P}_{\vartheta 0}:=\left\{x=(x_{kl})\in\Omega:\vartheta-\lim_{m,n\to\infty}\frac{1}{mn}\sum_{k,l=1}^{m,n}x_{kl}=0\right\}.
\end{eqnarray*}
Then, the sets $\mathcal{P}_{\infty}$, $\mathcal{P}_{\vartheta}$ and $\mathcal{P}_{\vartheta 0}$ can be defined by 4D matrix $T=(t_{mnkl})$ domain in the sequence space $\mathcal{BS}$, $\mathcal{CS}_{\vartheta}$ and $\mathcal{CS}_{\vartheta 0}$ where
\begin{eqnarray}
t_{mnkl}:=\left\{
\begin{array}{ccl}
\frac{1}{mn}&, & 1\leq k\leq m,~ 1\leq l\leq n, \\
0&, & \textrm{else}
\end{array}\right.
\end{eqnarray} 
Consider $\vartheta=\{p,p0,bp,bp0,r,r0\}$. Then, the spaces $\mathcal{P}_{\infty}=(\mathcal{BS})_{T}$, $\mathcal{P}_{\vartheta}=(\mathcal{CS}_{\vartheta})_{T}$ and $\mathcal{P}_{\vartheta 0}=(\mathcal{CS}_{\vartheta 0})_{T}$ are Banach spaces endowed with their natural norms expressed as follows: 
\begin{eqnarray}
\|x\|_{\mathcal{P}_{\infty}}=\sup_{m,n\in\mathbb{N}}\frac{1}{mn}\left|\sum_{k,l=1}^{m,n}x_{kl}\right|,
\end{eqnarray}
respectively. Since $\mathcal{C}_{p}\setminus\mathcal{M}_{u}$ is not empty, we can easily see that the following inclusions
\begin{eqnarray}\label{T domain of P}
\mathcal{P}_{\vartheta}\subset\mathcal{P}_{\infty}, ~~\mathcal{BS}\subset (\mathcal{BS})_{T},~~\mathcal{CS}_{\vartheta}\subset (\mathcal{CS}_{\vartheta})_{T}
\end{eqnarray}
strictly hold, where $\vartheta=\{bp,bp0,r,r0\}$. 

Now we give a literature review focusing on the Hahn sequence space $h$
from a single sequence perspective, before starting to the investigation of the Hahn double sequence space.

Hahn \cite{HH} defined the space $h$ as follows:
\begin{eqnarray}\label{hahn space}
h=\left\{x=(x_k)\in\omega:\sum_{k=1}^{\infty}k|\Delta x_k|<\infty\right\}\cap c_0,
\end{eqnarray}
where $\Delta$ represents the forward difference operator, i.e., $\Delta x_k=x_k-x_{k+1}$ for all $k=(1,2,3,...)$. Here, $\omega$ denotes the set of sequences with complex values, and $c_0$ refers the set of sequences converging to $0$. Hahn \cite{HH} proved that $h$ is a $BK$ when equipped with the norm
\begin{eqnarray}\label{norm of Hahn space}
\|x\|=\sum_{k=1}^{\infty}k|\Delta x_k|+\sup_{k\in\mathbb{N}}|x_k|\textit{ for all }x=(x_k)_{k=1}^{\infty}\in h.
\end{eqnarray}
Rao \cite{WCR} extended the study of $h$ and established that the Banach space $h$ has $AK$ property with the norm

\begin{eqnarray}\label{Rao norm of Hahn space}
	\|x\|=\sum_{k=1}^{\infty}k|\Delta x_k|<\infty \textit{ for all }x=(x_k)_{k=1}^{\infty}\in h.
\end{eqnarray}
Goes \cite{GGSG} investigated the generalized Hahn sequence space $h_d$ by considering an arbitrary negative or positive but not zero sequences $d=(d_k)_{k=1}^{\infty}$ for all $k\in\mathbb{N}$. The space $h_d$ is defined as follows:
\begin{eqnarray}\label{Generalized hahn space}
	h_d=\left\{x=(x_k)\in\omega:\sum_{k=1}^{\infty}|d_k||\Delta x_k|<\infty\right\}\cap c_0.
\end{eqnarray}
Most recently, Malkowsky et al \cite{TMR2} further extended the concept of the Hahn sequence space by studying $h_d$ with unbounded sequence $d=(d_k)_{k=1}^{\infty}$, where $d_k<d_{k+1}$ and $d_k>0$ for all $k\in\mathbb{N}$. Then Tu\v{g} et al \cite{Orhan aaa} studied the generalized $p-$Hahn sequence space $h_d^p$, where $1<p<\infty$.

\section{The Hahn double sequence space $\mathcal{H}_{\vartheta}$}
In the current section,  we begin by introducing the Hahn double sequence space $\mathcal{H}_{\vartheta}$. Then we study some topological structures exhibited by $\mathcal{H}_{\vartheta}$. Additionally, we establish several strict inclusion relations among the known double sequence spaces within the space $\mathcal{H}_{\vartheta}$.

We define Hahn double sequence space, denoted as $\mathcal{H}_{\vartheta}$, in the following manner:
\begin{eqnarray*}
\mathcal{H}_{\vartheta}:=\left\{x=(x_{kl})\in\Omega:\sum_{k,l=1}^{\infty}\left|kl\Delta x_{kl}\right|<\infty\right\}\cap\mathcal{C}_{\vartheta 0},
\end{eqnarray*}
where $\Delta$ represents the 4D forward difference operator, i.e., $\Delta x_{jk}=x_{jk}-x_{j+1,k}-x_{j,k+1}+x_{j+1,k+1}$, and $\vartheta\in\{p,bp,r\}$.

\begin{thm}\label{thm1}
\begin{enumerate}
\item[(i)] The spaces $\mathcal{H}_{bp}$ and $\mathcal{H}_{r}$ are complete normed spaces due to the norm
\begin{eqnarray}\label{Hahn norm}
	\|x\|_{\mathcal{H}_{\vartheta}}=\sum_{k,l=1}^{\infty}\left|kl\Delta x_{kl}\right|<\infty.
\end{eqnarray}
\item[(ii)] The space $\mathcal{H}_{p}$ is a complete semi-normed space, where the semi-norm is given by
\begin{eqnarray}
	\|x\|_{\mathcal{H}_{p}}=\sum_{k,l=1}^{\infty}\left|kl\Delta x_{kl}\right|+\lim_{N\to\infty}\left(\sup_{k,l\geq N}|x_{kl}|\right)
\end{eqnarray}
\end{enumerate}

\end{thm}

\begin{proof}
$(i)$ Let $\vartheta=\{bp,r\}$. Showing that the norm defined in equation (\ref{Hahn norm}) for $\mathcal{H}_{\vartheta}$ satisfies the norm axioms is straightforward. Therefore, we only investigate the completeness of the space $\mathcal{H}_{\vartheta}$. 

Define a Cauchy sequence $x^{[i]}=(x_{kl}^{[i]})\in \mathcal{H}_{\vartheta}$. Then it gives us
\begin{eqnarray}
\nonumber\sum_{k,l=0}^{\infty}\left|kl\Delta x_{kl}^{[i]}\right|<\infty,
\end{eqnarray}
\begin{eqnarray}
\nonumber\vartheta-\lim_{k,l\to\infty}x_{kl}^{[i]}=0,~\textit{ for every fixed } i\in\mathbb{N},
\end{eqnarray}
and $\forall\epsilon>0$, $\exists N(\epsilon)\in\mathbb{N}$ in which
\begin{eqnarray}
\label{eq1.11}\|x^{[i]}-x^{[j]}\|_{\mathcal{H}_{\vartheta}}<\frac{\epsilon}{2}
\end{eqnarray}
for every $i,j\geq N(\epsilon)$. For the fixed natural numbers $k,l$, and for all $q\geq k, p\geq l$ and all $i,j\geq N(\epsilon).$ Thus, one obtains
\begin{eqnarray}
\label{eq1.111}\sum_{m,n=k,l}^{q,p}(\Delta x_{mn}^{[i]}-\Delta x_{mn}^{[j]})\leq\sum_{m,n=k,l}^{q,p}\left|kl\Delta (x_{mn}^{[i]}- x_{mn}^{[j]})\right|\leq\|x^{[i]}-x^{[j]}\|_{\mathcal{H}_{\vartheta}}.
\end{eqnarray}
Thus, it can be deduced from (\ref{eq1.11}) and (\ref{eq1.111}) that
\begin{eqnarray*}
\left|x_{kl}^{[i]}-x_{kl}^{[j]}\right|&=&\left|\sum_{m,n=k,l}^{q,p}(\Delta x_{mn}^{[i]}-\Delta x_{mn}^{[j]})+(x_{q+1,p}^{[i]}-x_{q+1,p}^{[j]})\right.\\
&&\left.+(x_{q,p+1}^{[i]}-x_{q,p+1}^{[j]})-(x_{q+1,p+1}^{[i]}-x_{q+1,p+1}^{[j]})\right|\\
&<&\frac{\epsilon}{2}+\left|x_{q+1,p}^{[i]}\right|+\left|x_{q+1,p}^{[j]}\right|+\left|x_{q,p+1}^{[i]}\right|+\left|x_{q,p+1}^{[j]}\right|+\left|x_{q+1,p+1}^{[i]}\right|+\left|x_{q+1,p+1}^{[j]}\right|.
\end{eqnarray*}
Now, letting limit as $q,p\to\infty,$ according to the $\vartheta-$convergence, and since $x^{[i]},x^{[j]}\in\mathcal{C}_{\vartheta 0}$ for all $i,j\in\mathbb{N}$, one obtains
\begin{eqnarray*}
\left|x_{kl}^{[i]}-x_{kl}^{[j]}\right|<\epsilon
\end{eqnarray*}
for all $i,j\geq N(\epsilon)$. Consequently,  $x^{[i]}=(x_{kl}^{[i]})_{k,l=0}^{\infty}$ is a Cauchy sequence in $\mathbb{C}$ for every $k,l\in\mathbb{N}$. As $\mathbb{C}$ is complete, $x^{[i]}=(x_{kl}^{[i]})$ converges, i.e.,
\begin{eqnarray}\label{eq1.12}
\vartheta-\lim_{i\to\infty}x_{kl}^{[i]}=x_{kl},
\end{eqnarray}
where $x_{kl}\in\mathbb{C}$ for all natural numbers $k,l$. Let $q,p\in\mathbb{N}$ be given. Then we obtain by (\ref{eq1.11}) and (\ref{eq1.111}),
\begin{eqnarray}
\label{eq1.112}\sum_{k,l=1}^{q,p}\left|kl\Delta(x_{kl}^{[i]}-x_{kl}^{[j]})\right|\leq\|x^{[i]}-x^{[j]}\|_{\mathcal{H}_{\vartheta}}<\frac{\epsilon}{2}
\end{eqnarray}
for all $i,j\in\mathbb{N}$. Thus, the following inequality is obtained by letting limit on (\ref{eq1.112}) as $j\to\infty$ and by (\ref{eq1.12}),
\begin{eqnarray*}
	\sum_{k,l=1}^{q,p}\left|kl\Delta(x_{kl}^{[i]}-x_{kl})\right|\leq\frac{\epsilon}{2}.
\end{eqnarray*}
Since $q,p\in\mathbb{N}$ were arbitrary, we can write the following

\begin{eqnarray}\label{eq1.13}
\|x^{[i]}-x\|_{\mathcal{H}_{\vartheta}}<\epsilon
\end{eqnarray}
for all $i\geq N(\epsilon)$.

To show that $x\in\mathcal{C}_{\vartheta 0}$, we obtain the following inequality by using (\ref{eq1.12}) and by letting $\vartheta-$limit as $k,l\to\infty$ for every $i\geq N(\epsilon)$ that
\begin{eqnarray*}
0\leq |x_{kl}|\leq\left|x_{kl}-x_{kl}^{[i]}\right|+\left|x_{kl}^{[i]}\right|\to 0.
\end{eqnarray*}
Therefore, $x\in\mathcal{C}_{\vartheta 0}$. By utilizing (\ref{eq1.13}), it follows for every $n\geq N(\epsilon)$ that

\begin{eqnarray*}
\|x\|_{\mathcal{H}_{\vartheta}}\leq\left\|x-x^{[i]}\right\|_{\mathcal{H}_{\vartheta}}+\left\|x^{[i]}\right\|_{\mathcal{H}_{\vartheta}}<\infty.
\end{eqnarray*}
Thus, we conclude the proof by $x\in\mathcal{H}_{\vartheta}.$

$(ii)$ It closely resembles the proof of Theorem \ref{thm1}$(i)$. Hence, we choose to omit it.

\end{proof}

\begin{thm}\label{thm3}
The space $\mathcal{H}_{\vartheta}$, where $\vartheta\in\{p,bp,r\}$, is an $AK(\vartheta)-$space.
\end{thm}

\begin{proof}
Suppose that $\vartheta\in\{p,bp,r\},$ and that $x\in \mathcal{H}_{\vartheta}$. Then, it follows that
\begin{eqnarray*}
\|x-x^{[m-1,n-1]}\|_{\mathcal{H}_{\vartheta}}&=&\sum_{k,l=m,1}^{\infty,n-1}|kl\Delta x_{kl}|+\sum_{k,l=1,n}^{m-1,\infty}|kl\Delta x_{kl}|+\sum_{k,l=m,n}^{\infty}|kl\Delta x_{kl}|\\
&=&\sum_{k=m}^{\infty}\left(\sum_{l=1}^{n-1}|kl\Delta x_{kl}|\right)+\sum_{l=n}^{\infty}\left(\sum_{k=1}^{m-1}|kl\Delta x_{kl}|\right)+\sum_{k,l=m,n}^{\infty}|kl\Delta x_{kl}|
\end{eqnarray*}
for all $m,n\in\mathbb{N}$. Then,
\begin{eqnarray}\label{eq2.22}
\vartheta-\lim_{m,n\to\infty}\sum_{k,l=m,n}^{\infty}|kl\Delta x_{kl}|=0
\end{eqnarray}
clearly holds. Let us define the sequences $u=(u_{ml})$ and $v=(v_{kn})$  as
\begin{eqnarray*}
&&u_{ml}=\sum_{k=1}^{m-1}|kl\Delta x_{kl}|,\\
&&v_{kn}=\sum_{l=1}^{n-1}|kl\Delta x_{kl}|
\end{eqnarray*} 
for all $m,n,k,l\in\mathbb{N}$, respectively. Then, one obtains
\begin{eqnarray}
	\label{eq2.23}&&\sum_{l=n}^{\infty}\left(\sum_{k=1}^{m-1}|kl\Delta x_{kl}|\right)=\sum_{l=n}^{\infty}u_{ml},\\
	\label{eq2.24}&&\sum_{k=m}^{\infty}\left(\sum_{l=1}^{n-1}|kl\Delta x_{kl}|\right)=\sum_{k=m}^{\infty}v_{kn}.
\end{eqnarray}
Since the series $\sum_{k,l=1}^{\infty}|kl\Delta x_{kl}|$ converges, then by (\ref{eq2.23}) and (\ref{eq2.24}) the series $\sum_{l=n}^{\infty}u_{ml}$ and $\sum_{k=m}^{\infty}v_{kn}$ also converge. Therefore, the general terms of the series approach to zero as $m,l\to\infty$ and $k,n\to\infty$, respectively. Hence, we conclude the proof by using (\ref{eq2.22}) that
\begin{eqnarray*}
\vartheta-\lim_{m,n\to\infty}\|x-x^{[m-1,n-1]}\|_{\mathcal{H}_{\vartheta}}=0,
\end{eqnarray*}
as desired.
\end{proof}

\begin{thm}\label{thm2.1}
Each of the given assertions is strictly true.
\begin{enumerate}
\item[(i)] $\mathcal{H}_{\vartheta}\subset \mathcal{BV}$.
\item[(ii)] $\mathcal{H}_{\vartheta}\subset \mathcal{BV}_{\vartheta 0}$.
\end{enumerate}
\end{thm}

\begin{proof}
(i) Assume that $x=(x_{kl})$ is in the space $\mathcal{H}_{\vartheta}$. Then the sequence $x=(x_{kl})$ will be in the space $\mathcal{BV}$ since the following inequality

\begin{eqnarray}\label{eq2.2}
	\sum_{k,l=1}^{\infty}|\Delta x_{kl}|\leq \sum_{k,l=1}^{\infty}|kl\Delta x_{kl}|<\infty
\end{eqnarray}
 holds for every $k,l=1,2,3,...$. Therefore, the inclusion $\mathcal{H}_{\vartheta}\subset \mathcal{BV}$ holds. Now, we must show that the set $\mathcal{BV}\setminus\mathcal{H}_{\vartheta}$ is not empty. Consider $x=e$, then it is clear that $x\in\mathcal{BV}\setminus\mathcal{H}_{\vartheta}$.

(ii) First we show that the inclusion $\mathcal{H}_{\vartheta}\subset \mathcal{BV}_{\vartheta 0}$ holds. To do that, let us take a sequence $x=(x_{kl})\in\mathcal{H}_{\vartheta}$, that is, $x\in\mathcal{C}_{\vartheta 0}$ and $\sum_{k,l=1}^{\infty}|kl\Delta x_{kl}|<\infty$.  Since we can write the inequality (\ref{eq2.2}), then we can clearly see that $x\in\mathcal{BV}$. Thus, $x\in\mathcal{BV}_{\vartheta 0}$.

Now, to prove that the inclusion $\mathcal{H}_{\vartheta}\subset \mathcal{BV}_{\vartheta 0}$ strictly hold, we must show that the set $\mathcal{BV}_{\vartheta 0}\setminus\mathcal{H}_{\vartheta}$ is not empty. Consider $x=(x_{kl})$ as given in the following expression:
\begin{eqnarray}\label{ex2.1}
x_{kl}=
\begin{cases}
	\dfrac{1}{l} & (k=1)\\[1ex]
	0            & (k\ge 2) 
\end{cases}\quad(l=1,2,\dots).
\end{eqnarray}
Clearly $x\in\mathcal{C}_{\vartheta 0}$. Therefore, the $\Delta$ transform of the sequence $x$ is calculated as
\begin{eqnarray*}
\Delta x_{kl}=
\begin{cases}
	\dfrac{1}{l(l+1)} & (k=1)\\[2ex]
	0                           & (k\ge 2) 
\end{cases}\quad (l=1,2,\dots).
\end{eqnarray*}
Hence, one obtains
\begin{eqnarray*}
\sum\limits_{k,l=1}^{\infty}|\Delta x_{kl}|=
\sum\limits_{l=1}^{\infty}\left|\dfrac{1}{l(l+1)}\right|=1<\infty.
\end{eqnarray*}
But, it is apparent that
\begin{eqnarray*}
\sum\limits_{k,l=1}^{\infty}kl|\Delta x_{kl}|=
\sum\limits_{l=1}^{\infty}l\left|\dfrac{1}{l(l+1)}\right|=
\sum\limits_{l=1}^{\infty}\dfrac{1}{l+1}=\infty.
\end{eqnarray*}
Therefore, $x\in \mathcal{BV}$, but $x\notin \mathcal{H}_{\vartheta}$. It completes the proof.
\end{proof}

\begin{thm}
Neither the space $\mathcal{H}_{\vartheta}$ nor the space $\mathcal{L}_u$ include each other.
\end{thm}

\begin{proof}
We provide two counter-examples to establish the theorem.
\begin{enumerate}
\item[(i)] Choose a sequence $x=(x_{kl}),$ defined for all $k,l\in\mathbb{N}$ as follows:
\begin{eqnarray*}
	x_{kl}=
	\begin{cases}
		\dfrac{1}{kl} & (k=l)\\[1ex]
		0            & othervise .
	\end{cases}
\end{eqnarray*}
 It is clear that $x\in\mathcal{C}_{\vartheta 0}$. Moreover, we obtain that
\begin{eqnarray*}
	\Delta x_{kl}=
	\begin{cases}
		\dfrac{3kl+k+l+1}{kl(k+1)(l+1)} & (k=1, \forall l\in\mathbb{N})\\[2ex]
		2\left(\dfrac{2kl+k+l+1}{kl(k+1)(l+1)}\right) & (k\geq 2, \forall l\in\mathbb{N} )
	\end{cases}
\end{eqnarray*}
for all $k,l\in\mathbb{N}$. Therefore,
\begin{eqnarray*}
	\sum\limits_{k,l=1}^{\infty}|x_{kl}|=
	\sum\limits_{k=l=1}^{\infty}\left|\dfrac{1}{kl}\right|=\sum\limits_{r=1}^{\infty}\dfrac{1}{r^2}=\frac{\pi^2}{6}<\infty,
\end{eqnarray*}
but
\begin{eqnarray*}
	\sum\limits_{k,l=1}^{\infty}kl|\Delta x_{kl}|&=&\sum\limits_{l=1}^{\infty}l\left|\dfrac{2l+1}{l(l+1)}\right|+\sum\limits_{k=2}^{\infty}2k\left|\dfrac{3k+2}{2k(k+1)}\right|+\sum\limits_{k=2}^{\infty}\sum\limits_{l=2}^{\infty}2kl\left|\dfrac{2kl+k+l+1}{kl(k+1)(l+1)}\right|=\infty.
\end{eqnarray*}
Thus, the set $\mathcal{L}_u\setminus\mathcal{H}_{\vartheta}$ is non-empty.

\item[(ii)] Take $x=(x_{kl})$ defined by the following expression: 
\begin{eqnarray*}
	x_{kl}=
	\begin{cases}
		1 & k=1, \forall l\in\mathbb{N}\\[1ex]
		0            & othervise 
	\end{cases}
\end{eqnarray*}
for all $k,l\in\mathbb{N}$. Clearly $x\in\mathcal{C}_{\vartheta 0}$. Moreover, we obtain that

\begin{eqnarray*}
	\sum\limits_{k,l=1}^{\infty}|x_{kl}|=
	\sum\limits_{l=1}^{\infty}\left|x_{1l}\right|=\sum\limits_{k,l=1}^{\infty}1=\infty.
\end{eqnarray*}
Thus, $x\notin\mathcal{L}_u$. However,

\begin{eqnarray*}
	\sum\limits_{k,l=1}^{\infty}kl|\Delta x_{kl}|=\sum\limits_{l=1}^{\infty}l\left|x_{1l}-x_{2l}-x_{1,l+1}+x_{2,l+1}\right|=\sum_{l=1}^{\infty}l|x_{1l}-x_{1,l+1}|=0<\infty,
\end{eqnarray*}
which shows $x\in\mathcal{H}_{\vartheta}$. Therefore, the set $\mathcal{H}_{\vartheta}\setminus\mathcal{L}_u$ is not empty, too. 

\end{enumerate}
 
\end{proof}

The following Corollary is immediate from the relations $\sum_{k,l=1}^{\infty}|kl\Delta x_{kl}|\leq4 \sum_{k,l=1}^{\infty}|klx_{kl}|$ and $(x_{kl})=\left\{\frac{1}{(kl)^2}\right\}\in \mathcal{H}_{\vartheta}\setminus\int\mathcal{L}_u$

\begin{cor}
The inclusion $\int\mathcal{L}_u\subset\mathcal{H}_{\vartheta}$ strictly holds.
\end{cor}

\section{The Dual Spaces of $\mathcal{H}_{\vartheta}$}

Within this section, we undertake the computation of the algebraic dual, as well as $\alpha-$, $\beta(bp)-$, and $\gamma-$ dual of the space $\mathcal{H}_{\vartheta}$. Before starting to the section's main content, we perform some elementary calculation and establish the following Lemmas, which are key in obtaining our mai results.

Let us take $kl\Delta x_{kl}=y_{kl}$, where $y\in\mathcal{L}_u$, $x\in\mathcal{H}_{\vartheta}$ and $k,l\in\mathbb{N}$. Then, we write
\begin{eqnarray*}
	\frac{1}{kl}y_{kl}&=&\Delta x_{kl}\\
	\frac{1}{(k+1)l}y_{k+1,l}&=&\Delta x_{k+1,l}\\
	\frac{1}{k(l+1)}y_{k,l+1}&=&\Delta x_{k,l+1}\\
	\frac{1}{(k+1)(l+1)}y_{k+1,l+1}&=&\Delta x_{k+1,l+1}\\
	&\vdots&\\
	\frac{1}{mn}y_{mn}&=&\Delta x_{mn}.
\end{eqnarray*}
Adding both sides of the above equations from ranging from $(k,l)$ to $(m,n)$, we obtain
\begin{eqnarray}\label{eq3.8}
	\sum_{i,j=k,l}^{m,n}\frac{y_{ij}}{ij}=x_{kl}-x_{m+1,l}-x_{k,n+1}+x_{m+1,n+1}.
\end{eqnarray}
Thus, by letting $\vartheta-$limit over the equality in (\ref{eq3.8}) as $m,n\to\infty$, then we get
\begin{eqnarray}\label{eq3.9}
	\vartheta-\lim_{m,n\to\infty}\sum_{i,j=k,l}^{m,n}\frac{y_{ij}}{ij}=\sum_{i,j=k,l}^{\infty}\frac{y_{ij}}{ij}=x_{kl}
\end{eqnarray}
since $\vartheta-\lim_{m,n\to\infty}x_{mn}=0$. 

Now, we write the following Lemma with the specific case for $s=1$ from the results \cite[Theorem 4.1.(i), Theorem 4.2., Theoerem 4.3.(i)]{Medine}, respectively, as follows:

\begin{lem}\label{lem2.1}
Each of the following assertions holds true:
\begin{enumerate}
\item[(i)] A 4D matrix $A\in(\mathcal{L}_u,\mathcal{M}_{u})$ iff 
\begin{equation}
	\label{eq3.3}\sup_{m,n,k,l}|a_{mnkl}|<\infty.
\end{equation}
\item[(ii)] A 4D matrix $A\in(\mathcal{L}_u,\mathcal{L}_{u})$ iff  
\begin{equation}
	\label{eq3.5}\sup_{k,l\in\mathbb{N}}\sum_{m,n}|a_{mnkl}|<\infty.
\end{equation}
\item[(iii)] A 4D matrix $A\in(\mathcal{L}_u,\mathcal{C}_{bp})$ iff  the condition in (\ref{eq3.3}) holds true, and there exists $a_{kl}\in\mathbb{C}$ satisfying
\begin{equation}
	\label{eq3.4}bp-\lim_{m,n\to\infty}a_{mnkl}=a_{kl}.
\end{equation}
\end{enumerate}
\end{lem}

\begin{thm}
Assume $\vartheta$ to be any one of the space $bp$, $p$ or $r$. Then, $\{\mathcal{H}_{\vartheta}\}^{\alpha}=h_1$, where
\begin{eqnarray}
h_1=\left\{ a=(a_{mn})\in\Omega:\sup_{k,l\in\mathbb{N}}\frac{1}{kl}\sum_{m,n=0}^{\infty}\left|a_{mn}\right|<\infty \right\}.
\end{eqnarray}
\end{thm}

\begin{proof}
Let $a=(a_{mn})\in \Omega$ be given double sequence and let us consider the following equality
\begin{eqnarray}\label{eq 4.56}
a_{mn}x_{mn}=\sum_{k,l=m,n}^{\infty}\frac{a_{mn}}{kl}y_{kl}=(Dy)_{mn}, \textit{ for every }m,n\in\mathbb{N},
\end{eqnarray}
where $D=(d_{mnkl})$ is defined by
\begin{eqnarray}
	d_{mnkl}:=\left\{
	\begin{array}{ccl}
		\frac{a_{mn}}{kl}&, & k\geq m,~ l\geq n, \\
		0&, & \textrm{otherwise}
	\end{array}\right.
\end{eqnarray}	
for every natural number $m,n,k,l$. It follows from the equality in (\ref{eq 4.56}) with Lemma \ref{lem2.1}(ii) that $ax=(a_{mn}x_{mn})\in\mathcal{L}_u$ whenever $x=(x_{kl})\in\mathcal{H}_{\vartheta}$ iff  $Dy\in\mathcal{L}_u$ whenever $y=(y_{kl})\in\mathcal{L}_u$. This shows that $a\in\{\mathcal{H}_{\vartheta}\}^{\alpha}$ whenever $x=(x_{kl})\in\mathcal{H}_{\vartheta}$ iff  $D\in(\mathcal{L}_u,\mathcal{L}_u)$. That is to say that the condition stated as in (\ref{eq3.5}) is satisfied with $d_{mnkl}$ replacing $a_{mnkl}$. This proves that $\{\mathcal{H}_{\vartheta}\}^{\alpha}=h_1$. 
\end{proof}

Now we give the following Lemma for double sequences (cf. Wilansky \cite[Theorem 7.2.9]{WA} for single sequences).

\begin{lem}\label{lemma 2}
Let $E$ be a $BDK-$space with $AK(\vartheta)$, and $E^*$ be the space of linear continuous functional on $E$. Then, $E^*$ can be identified with $E^{\beta(\vartheta)}$ by the isomorphism 
\begin{eqnarray}\label{eq 4.7}
	T&:&E^*\to E^{\beta(\vartheta)}\\
	&&\nonumber f\to T(f)=\{f(e^{kl})\}_{k,l=1}^{\infty}.
\end{eqnarray}
Because, $E$ has $AK(\vartheta)$ iff  every $f\in E^*$ can be written in the form $f(a)=\sum_{k,l=1}^{\infty}a_{kl}b_{kl}$, where $a=(a_{kl})\in E$ and $b=(b_{kl})\in E^{\beta(\vartheta)}$.
\end{lem}

\begin{proof}
First we show that $T$ is a well-defined operator. Suppose $f \in E^*$, we must show that $T(f) \in E^{\beta(\vartheta)}$, that is, the sequence $T(f)$ is $\beta(\vartheta)$-convergent. Since $E$ is a $BDK$ space with $AK(\vartheta)$, we can say that for each $f \in E^*$, there exists a sequence $a \in E$ in which 
\begin{eqnarray*}
f(x) = \sum_{k,l=1}^{\infty} a_{kl}b_{kl}
\end{eqnarray*}
for some sequences $b \in E^{\beta(\vartheta)}$. Therefore, 
\begin{eqnarray*}
f(e^{kl}) = \sum_{k,l=1}^{\infty} a_{kl}b_{kl}.
\end{eqnarray*}

Now, let's show that $T(f)$ is $\vartheta$-convergent. Actually, this is clearly seen from the fact that each component of $f(e^{kl})$ is a convergent sequence in $E^{\beta(\vartheta)}$. Thus, $T(f) \in E^{\beta(\vartheta)}$.

Moreover, linearity of $T$ is clear. To show that $T$ is an isomorphism, we must show that it is bijective. Suppose $T(f) = T(g)$ for some $f, g \in E^*$. Then, for any $x \in E$, $T(f)(x) = T(g)(x)$, which implies that $f(x) = g(x)$. Since $x$ is arbitray, therefore, $T$ is injective. Now, let $h \in E^{\beta(\vartheta)}$. Since $E$ has $AK(\vartheta)$, for any $h \in E^{\beta(\vartheta)}$, there exists $f \in E^*$ in which $T(f) = h$. We have this because $E$ has $AK(\vartheta)$ iff every $f\in E^*$ can be written as $f(a) = \sum_{k,l=1}^{\infty} a_{kl}b_{kl}$, where $a \in E$ and $b \in E^{\beta(\vartheta)}$. Therefore, $T$ is surjective, too. Hence, it is now clearly seen that $T: E^* \to E^{\beta(\vartheta)}$ is a well-defined linear isomorphism.
\end{proof}

\begin{thm}
Let $\vartheta=\{bp,r\}$. The continuous dual $\mathcal{H}_{\vartheta}^*$ and $\beta(\vartheta)-$dual $\{\mathcal{H}_{\vartheta}\}^{\beta(\vartheta)}$ of the space $\mathcal{H}_{\vartheta}$ are norm isomorphic.
\end{thm}

\begin{proof}
Suppose that $\vartheta=\{bp,r\}$. It is established that $\mathcal{H}_{\vartheta}$ is a $BDK-$space and has $AK(\vartheta)$ by Theorem \ref{thm1} and Theorem \ref{thm3}, respectively. Moreover, the spaces $\mathcal{H}_{\vartheta}^*$ and $\mathcal{H}_{\vartheta}^{\beta(\vartheta)}$ are isomorphic by Remark \ref{lemma 2}. Therefore, it is sufficient to prove the fact that
\begin{eqnarray*}
\|f\|=\|a\|_{\mathcal{P}_{\infty}}
\end{eqnarray*}
for all $f\in \mathcal{H}_{\vartheta}^*$ represented by the sequence $a\in \mathcal{P}_{\infty}$ in which 
\begin{eqnarray*}
f(x)=\sum_{k,l=1}^{\infty}a_{ kl}x_{kl}\textit{ for all }x\in\mathcal{H}_{\vartheta}.
\end{eqnarray*}

Suppose that $f\in \mathcal{H}_{\vartheta}^*$ and that $a=(a_{kl})\in\mathcal{P}_{\infty}$ in which  $f(x)=\sum_{k,l=1}^{\infty}a_{ kl}x_{kl} \text{ for all }x\in\mathcal{H}_{\vartheta}$. Moreover, let $m,n\in\mathbb{N}$  and $x\in\mathcal{H}_{\vartheta}.$ Therefore, by utilizing Abel's double summability by parts, it follows that
\begin{eqnarray*}
f(x^{[m,n]})=\sum_{k,l=1}^{m,n}a_{ kl}x_{kl}&=&\sum_{k,l=1}^{m-1,n-1}\Delta_{11}^{kl}x_{kl}\sum_{i,j=1}^{k,l}a_{ij}+\sum_{l=1}^{n-1}\Delta_{01}^{kl}x_{ml}\sum_{i,j=1}^{m,l}a_{ij}\\
&&+\sum_{k=1}^{m-1}\Delta_{10}^{kl}x_{kn}\sum_{i,j=1}^{k,n}a_{ij}+x_{mn}\sum_{i,j=1}^{m,n}a_{ij}\\
&=&\sum_{k,l=1}^{m-1,n-1}kl\Delta_{11}^{kl}x_{kl}\left(\frac{1}{kl}\sum_{i,j=1}^{k,l}a_{ij}\right)+\sum_{l=1}^{n-1}ml\Delta_{01}^{kl}x_{ml}\left(\frac{1}{ml}\sum_{i,j=1}^{m,l}a_{ij}\right)\\
&&+\sum_{k=1}^{m-1}kn\Delta_{10}^{kl}x_{kn}\left(\frac{1}{kn}\sum_{i,j=1}^{k,n}a_{ij}\right)+mnx_{mn}\left(\frac{1}{mn}\sum_{i,j=1}^{m,n}a_{ij}\right),
\end{eqnarray*}
where
\begin{eqnarray*}
	\Delta_{10}^{kl}x_{kl}&=&x_{kl}-x_{k+1,l},\\
	\Delta_{01}^{kl}x_{kl}&=&x_{kl}-x_{k,l+1},\\
	\Delta_{11}^{kl}x_{kl}&=& \Delta_{10}^{kl}(\Delta_{01}^{kl}x_{kl})=\Delta_{01}^{kl}(\Delta_{10}^{kl}x_{kl})=\Delta x_{kl}
\end{eqnarray*}
It is known by Theorem \ref{thm3} that the space $\mathcal{H}_{\vartheta}$ has $AK(\vartheta)$, and given  $x\in\mathcal{H}_{\vartheta}$, we obtain that
\begin{eqnarray*}
0\leq|mnx_{mn}|&=&\sum_{k,l=m,n}^{\infty}kl|\Delta x_{kl}|\\
&\leq&\sum_{k,l=1}^{\infty}kl\left|\Delta (x_{kl}^{[m,n]}-x_{kl})\right|+\sum_{k,l=m,n}^{\infty}kl\left|\Delta x_{kl}\right|\\
&=&\left\|x-x^{[m,n]}\right\|_{\mathcal{H}_{\vartheta}}+\sum_{k,l=m,n}^{\infty}kl\left|\Delta x_{kl}\right|\to 0 \textit{ as }(m,n\to\infty).
\end{eqnarray*}
Hence, by using the continuity of $f$, we get that
\begin{eqnarray*}
|f(x)|=\lim_{m,n\to\infty}\left|f(x^{[m,n]})\right|\leq\left(\sum_{k,l=1}^{\infty}kl|\Delta x_{kl}|\right)\|a\|_{\mathcal{P}_{\infty}}=\|x\|_{\mathcal{H}_{\vartheta}}\|a\|_{\mathcal{P}_{\infty}}
\end{eqnarray*}
and hence
\begin{eqnarray}\label{isomorphism 1}
\|f\|\leq\|a\|_{\mathcal{P}_{\infty}}.
\end{eqnarray}

Now, assume that $m,n\in\mathbb{N}$, and that $x^{(m,n)}=\frac{e^{[m,n]}}{mn}$ which clearly belongs to the space $\mathcal{H}_{\vartheta}$ and 
\begin{eqnarray}\label{isomorphism 2}
\|x^{(m,n)}\|_{\mathcal{H}_{\vartheta}}=\sum_{k,l=1}^{\infty}kl\left|\Delta x_{kl}^{(m,n)}\right|=1.
\end{eqnarray}
Therefore,
\begin{eqnarray*}
\left|\sum_{k,l=1}^{\infty}a_{kl}x_{kl}^{(m,n)}\right|=\frac{1}{mn}\left|\sum_{k,l=1}^{m,n}a_{kl}\right|\leq \|f\|.
\end{eqnarray*}
Since $m,n\in\mathbb{N}$ were arbitrary and (\ref{isomorphism 2}) holds, we obtain that
\begin{eqnarray}\label{isomorphism 3}
\|a\|_{\mathcal{P}_{\infty}}\leq\|f\|.
\end{eqnarray}
Hence, we get the equality $\|a\|_{\mathcal{P}_{\infty}}=\|f\|$ by the results (\ref{isomorphism 1}) and (\ref{isomorphism 3}). The proof is thereby finished.
\end{proof}

\begin{thm}\label{betadual}
$\{\mathcal{H}_{\vartheta}\}^{\beta(bp)}=\mathcal{P}_{\infty}$.
\end{thm}

\begin{proof}
Suppose that $a=(a_{kl}) \in \Omega$ and $x=(x_{kl})\in\mathcal{H}_{\vartheta}$. Then, we have
\begin{eqnarray}\label{betadual 1}
\sum_{k,l=1}^{m,n}a_{kl}x_{kl}=\sum_{k,l=1}^{m,n}a_{kl}\left(\sum_{i,j=k,l}^{m,n}\frac{y_{ij}}{ij}\right)=\sum_{k,l=1}^{m,n}\left(\frac{1}{kl}\sum_{i,j=1}^{k,l}a_{ij}\right)y_{kl}=(By)_{mn}
\end{eqnarray}
for all positive integers $m,n$, where $B=(b_{mnkl})$ is defined by
\begin{eqnarray}
	b_{mnkl}:=\left\{
	\begin{array}{ccl}
		\sum_{i,j=1}^{k,l}\frac{a_{ij}}{kl}&, & k\geq m,~ l\geq n, \\
		0&, & \textrm{otherwise}
	\end{array}\right.
\end{eqnarray}	
for all natural numbers $m,n,k,l$. Then, by using (\ref{betadual 1}) and Lemma \ref{lem2.1}(iii), we obtain  that $ax\in\mathcal{CS}_{bp}$ whenever $x=(x_{kl})\in\mathcal{H}_{\vartheta}$ iff  $By\in\mathcal{C}_{bp}$ whenever $y\in\mathcal{L}_{u}$. Therefore, $a\in\{\mathcal{H}_{\vartheta}\}^{\beta(bp)}$ iff  $B\in(\mathcal{L}_{u},\mathcal{C}_{bp}).$ As a result the conditions in (\ref{eq3.3}) and (\ref{eq3.4}) hold with $b_{mnkl}$ instead of $a_{mnkl}$, that is,
\begin{eqnarray}
	\label{eq3.30}\sup_{m,n,k,l\in\mathbb{N}}\frac{1}{kl}\left|\sum_{i,j=1}^{k,l}a_{ij}\right|<\infty\\
		\label{eq3.40}bp-\lim_{k,l\to\infty}\frac{1}{kl}\sum_{i,j=1}^{k,l}a_{ij} \textit{ exists}.
\end{eqnarray}
This shows that $a=(a_{kl})\in\mathcal{P}_{\infty}$ and $a=(a_{kl})\in\mathcal{P}_{bp}$. Moreover, it is easy to see by (\ref{T domain of P}) that $\mathcal{P}_{bp}\subset\mathcal{P}_{\infty}$. Thus, we  conclude that $\{\mathcal{H}_{\vartheta}\}^{\beta(bp)}=\mathcal{P}_{\infty}$.
\end{proof}

\begin{thm}
$\{\mathcal{H}_{\vartheta}\}^{\gamma}=\mathcal{P}_{\infty}$. 
\end{thm}

\begin{proof}
Utilizing a comparable approach to that used in Theorem \ref{betadual}, the proof can be summarized as follows:

we have from the equality in (\ref{betadual 1}) and Lemma \ref{lem2.1}(i) that $ax\in\mathcal{BS}$ whenever $x=(x_{kl})\in\mathcal{H}_{\vartheta}$ iff  $By\in\mathcal{M}_{u}$ whenever $y\in\mathcal{L}_{u}$. Therefore, $a\in\{\mathcal{H}_{\vartheta}\}^{\gamma}$ iff  $B\in(\mathcal{L}_{u},\mathcal{M}_{u})$
in which  the condition stated in (\ref{eq3.3}) is satisfied with $b_{mnkl}$ replacing $a_{mnkl}$. Therefor, the condition in (\ref{eq3.30}) holds which shows that $a=(a_{kl})\in\mathcal{P}_{\infty}$. Thus, it gives the result that $\{\mathcal{H}_{\vartheta}\}^{\gamma}=\mathcal{P}_{\infty}$.
\end{proof}

\section{Matrix Mappings from and into the Space $\mathcal{H}_{\vartheta}$}

Within this section, we represent the matrix transformations that operate both from and into the sequence space $\mathcal{H}_{\vartheta}$. First, we characterize the matrix classes $(\mathcal{H}_{\vartheta},\lambda)$, where $\lambda=\{\mathcal{M}_u,\mathcal{L}_u, \mathcal{C}_{bp}\}$. Then we define determining set of $\mathcal{H}_{\vartheta}$ and then we characterize the matrix class $(\mathcal{H}_{\vartheta},\mathcal{H}_{\vartheta})$. Moreover, we determine essential conditions of the matrix classes $(\mathcal{H}_{\vartheta},\lambda)$, where $\lambda=\{\mathcal{BV}, \mathcal{BV}_{\vartheta 0}, \mathcal{CS}_{\vartheta},\mathcal{CS}_{\vartheta 0},\mathcal{BS}\}$ and $(\mu,\mathcal{H}_{\vartheta})$, where $\mu=\{\mathcal{L}_u, \mathcal{C}_{\vartheta 0}, \mathcal{C}_{\vartheta},\mathcal{M}_{u}\}$.

In the entirety of this text, we presume that $A=(a_{mnkl})$ is a 4D matrix over $\mathbb{C}$.

\begin{thm}\label{thm4-1}
$A=(a_{mnkl})\in(\mathcal{H}_{\vartheta},\lambda)$ iff 
\begin{eqnarray}
\label{thm4.1i}A_{mn}\in\mathcal{H}_{\vartheta}^{\beta(\vartheta)},\\
\label{thm4.1ii}E\in(\mathcal{L}_u,\lambda),
\end{eqnarray}
where the 4D matrix $E=(e_{mnkl})$ is given by
\begin{equation}\label{matrix4-1}
		e_{mnkl}:=\left\{
		\begin{array}{ccl}
			\frac{1}{kl}\sum_{i,j=1}^{k,l}a_{mnij}&, & k\geq m,~ l\geq n, \\
			0&, & \textrm{otherwise}
		\end{array}\right.
\end{equation}
\end{thm}

\begin{proof}
Assume that $A=(a_{mnkl})\in(\mathcal{H}_{\vartheta},\lambda)$. Consequently, $A_{mn}(x)$ exists for every $x=(x_{kl})\in\mathcal{H}_{\vartheta}$ and contains in $\lambda$ for every positive integer $m$ and $n$. The necessity of the condition stated in (\ref{thm4.1i}) is clear by considering the entries of the sequence $x=(x_{kl})$ as $x_{kl}=e^{kl}$ for all $k,l\in\mathbb{N}$.

Assume that $x=(x_{kl})\in\mathcal{H}_{\vartheta}.$ Then, we deduce the subsequent equation:
\begin{eqnarray}\label{eq4-1}
\sum_{k,l=1}^{m,n}a_{mnkl}x_{kl}=\sum_{k,l=1}^{m,n}a_{mnkl}\left(\sum_{i,j=k,l}^{m,n}\frac{y_{ij}}{ij}\right)=\sum_{k,l=1}^{m,n}\left(\sum_{i,j=1}^{k,l}\frac{a_{mnij}}{kl}\right)y_{kl}=(Ey)_{mn},
\end{eqnarray}
where $E=(e_{mnkl})$ is given as in (\ref{matrix4-1}). By letting limit as $m,n\to\infty$, we deduce the fact that $Ax=Ey$. In this case, we say that for every $x=(x_{kl})\in\mathcal{H}_{\vartheta}$, $Ax\in \lambda$ exists for all natural numbers $m,n$ iff  for every $y=(y_{kl})\in\mathcal{L}_{u}$, $Ey\in\lambda$ exists for all $m,n\in\mathbb{N}.$ That is to say that $A=(a_{mnkl})\in(\mathcal{H}_{\vartheta},\lambda)$ iff  $E\in(\mathcal{L}_u,\lambda)$, which indicates the requirement for (\ref{thm4.1ii}). 

Conversely, we assume that conditions in (\ref{thm4.1i}) and (\ref{thm4.1ii}) are satisfied and the relation (\ref{matrix4-1}) between the 4D infinite matrices $A=(a_{mnkl})$ and $E=(e_{mnkl})$ be given. Since the condition (\ref{thm4.1i}) holds, $Ax$ exists for every $x\in\mathcal{H}_{\vartheta}$. Moreover, the equation given in (\ref{eq4-1}) and the condition in (\ref{thm4.1ii}) show the existence of $Ax\in \lambda$. Therefore, $A=(a_{mnkl})\in(\mathcal{H}_{\vartheta},\lambda)$.
\end{proof}

\begin{cor}
Each of the following statements are valid:
\begin{enumerate}
\item[(i)] $A=(a_{mnkl})\in(\mathcal{H}_{\vartheta},\mathcal{M}_{u})$ iff  (\ref{thm4.1i}) holds, and the condition stated in (\ref{eq3.3}) is satisfied with $e_{mnkl}$ replacing $a_{mnkl}$.
\item[(ii)] $A=(a_{mnkl})\in(\mathcal{H}_{\vartheta},\mathcal{L}_{u})$ iff  (\ref{thm4.1i}) holds, and the condition stated in (\ref{eq3.5}) is satisfied with $e_{mnkl}$ replacing $a_{mnkl}$.
\item[(ii)] $A=(a_{mnkl})\in(\mathcal{H}_{\vartheta},\mathcal{C}_{bp})$ iff  (\ref{thm4.1i}) holds, and the conditions stated in (\ref{eq3.3}) and (\ref{eq3.4}) are satisfied with $e_{mnkl}$ replacing $a_{mnkl}$.
\end{enumerate}
\end{cor}

\begin{defin}\cite[Definition 7.4.2]{WA}
	Let $X$ be a $BK$ space. A subset $E$ of the set $\phi$ called a determining set for $X$ if $D(X)=\overline{B}_X\cap\phi$ is the absolutely convex hull of $E$.
\end{defin}

\begin{thm}\label{thm4.4}
Consider the set $E$ consisting all sequences in $\phi$, each having non-zero elements given by
	\begin{eqnarray*}
		\begin{pmatrix}
			0 & 0& \dots & 0 & \dots  & 0 \\
			0 & 0& \dots & 0 & \dots  & 0 \\
			\vdots & \vdots & \ddots & \vdots & \ddots & \vdots \\
			0 & 0& \dots & \frac{1}{mn} & \dots  & 0 \\
			0 & 0& \dots & 0 & \dots  & 0 \\
			\vdots & \vdots & \ddots & \vdots & \ddots & \vdots 
		\end{pmatrix}
	\end{eqnarray*}
where $\frac{1}{mn}$ occupies the $(m,n)^{th}-$place and all other entries are zero. Then it follows that $E$ is a determining set for $\mathcal{H}_{\vartheta}$.
\end{thm}

\begin{proof}
	Let $D(\mathcal{H}_{\vartheta})=\left\{x=(x_{kl})\in\phi:\|x\|_{\mathcal{H}_{\vartheta}}\leq1\right\}$, where $\phi$ represents the space of all finite double sequences, and $K[E]$ be the absolutely convex hull of $E$. We must prove that $D(\mathcal{H}_{\vartheta})=K[E]$ by showing that the inclusions 
	\begin{eqnarray}
		\label{konvex hull 1}D(\mathcal{H}_{\vartheta})\subset K[E],\\
		\label{konvex hull 2}K[E]\subset D(\mathcal{H}_{\vartheta})
	\end{eqnarray}
	hold. First we show the inclusion (\ref{konvex hull 1}). Suppose that $x\in D(\mathcal{H}_{\vartheta})$. Since $x\in\phi$, we obtain, for any positive inter $m$ and $n$ that 
	\begin{eqnarray*}
		x=x^{[m,n]}=\sum_{k,l=1}^{m,n}z_{kl}d_{kl} 
	\end{eqnarray*}
	holds, where $d_{kl}=\frac{1}{kl}e^{kl}$  and $z_{kl}=kl\Delta x_{kl}$ for any positive integer $k$ and $l$. This can inferred from that facts that $x_{m+1,l}=x_{k,n+1}=x_{m+1,n+1}=0$ for each $k,l\in\mathbb{N}$ and 
	\begin{eqnarray*}
		\sum_{k,l=1}^{m,n}z_{kl}d_{kl}&=&\sum_{k,l=1}^{m,n}kl\Delta x_{kl}d_{kl}\\
		&=&\sum_{k,l=1}^{m,n}kl x_{kl}d_{kl}-\sum_{k,l=1}^{m-1,n}kl x_{k+1,l}d_{kl}-\sum_{k,l=1}^{m,n-1}kl x_{k,l+1}d_{kl}+\sum_{k,l=1}^{m-1,n-1}kl x_{k+1,l+1}d_{kl}\\
		&=&x_{11}d_{11}+\sum_{l=2}^{n}lx_{1,l}d_{1,l}+\sum_{k=2}^{m}kx_{k,1}d_{k,1}+\sum_{k,l=2}^{m,n}klx_{kl}d_{kl}\\
		&&-\left(\sum_{k=2}^{m}(k-1)x_{k,1}d_{k-1,1}+\sum_{k,l=2}^{m,n}(k-1)lx_{kl}d_{k-1,l} \right)\\
		&&-\left(\sum_{l=2}^{n}(l-1)x_{1,l}d_{1,l-1}+\sum_{k,l=2}^{m,n}k(l-1)x_{kl}d_{k,l-1} \right)+\sum_{k,l=2}^{m,n}(k-1)(l-1)x_{kl}d_{k-1,l-1}\\
		&=&x_{11}e^{11}+\sum_{k=2}^{m}x_{k,1}(kd_{k,1}-(k-1)d_{k-1,1})+\sum_{l=2}^{n}x_{1,l}(ld_{1,l}-(l-1d_{1,l-1}))\\
		&&+\sum_{k,l=2}^{m,n}x_{kl}(kld_{kl}-(k-1)ld_{k-1,l}-k(l-1)d_{k,l-1}+(k-1)(l-1)d_{k-1,l-1})\\
		&=&x_{11}e^{1,1}+\sum_{k=2}^{m}x_{k,1}e^{k,1}+\sum_{l=2}^{n}x_{1,l}e^{1,l}+\sum_{k,l=2}^{m,n}x_{kl}e^{kl}\\
		&=&x^{[m,n]}.
	\end{eqnarray*}
	Moreover, it is observed that
	\begin{eqnarray*}
		\sum_{k,l=1}^{m,n}|z_{kl}|=\sum_{k,l=1}^{m,n}kl|\Delta x_{kl}|\leq \|x\|_{\mathcal{H}_{\vartheta}}\leq1.
	\end{eqnarray*}
	Therefore, $x\in K[E]$, which proves (\ref{konvex hull 1}).

	Now, suppose that $x\in K[E]$. Then there are natural numbers $m,n$ such that
	\begin{eqnarray*}
		x=\sum_{k,l=1}^{m,n}z_{kl}d_{kl} \textit{ with } \sum_{k,l=1}^{m,n}|z_{kl}|\leq 1.
	\end{eqnarray*}
	Since $x_{m+1,l}=x_{k,n+1}=x_{m+1,n+1}=0$ for each $k,l\in\mathbb{N}$, that is, $x_{kl}=0$ for all $k\geq m+1$ or $l\geq n+1$ or both, by using the relation in (\ref{eq3.8}) and for all $1\leq k\leq m$ and $1\leq l\leq n$ we obtain that
	\begin{eqnarray*}
		\sum_{i,j=k,l}^{m,n}\frac{z_{ij}}{ij}=\sum_{i,j=k,l}^{m,n}\Delta x_{ij}=x_{kl}-x_{m+1,l}-x_{k,n+1}+x_{m+1,n+1}=x_{kl}.
	\end{eqnarray*}
	Furthermore, since $x_{m+1,l}=x_{k,n+1}=x_{m+1,n+1}=0$ for each $k,l\in\mathbb{N},$ we have that 
	\begin{eqnarray*}
		\|x\|_{\mathcal{H}_{\vartheta}}=\sum_{k,l=1}^{\infty}kl|\Delta x_{kl}|=\sum_{k,l=1}^{m,n}|z_{kl}|\leq 1.
	\end{eqnarray*}
	Hence, $x\in  D(\mathcal{H}_{\vartheta})$. Thus, the inclusion (\ref{konvex hull 2}) is established. The result immediately follows from both the inclusions (\ref{konvex hull 1}) and (\ref{konvex hull 2}).
\end{proof}

Now we give the following Lemma for double sequence spaces with respect to the notation in \cite[Theorem 8.3.4]{WA} as follows:

\begin{lem}\label{Lemma 4.5}
	Let $X$ be a $BDK$ space with $AK(\vartheta)$, $E$ be a determining set for $X$, $Y$ be an $FDK$ space, and $A=(a_{mnkl})$ is a 4D matrix. Suppose that either $X$ has $AK(\vartheta)$ or $A$ is row finite. Then $A\in(X,Y)$ iff 
	\begin{enumerate}
		\item[(i)] The columns of $A$ belongs to $Y$, that is, $A^{kl}=(a_{mnkl})_{m,n=1}^{\infty}\in Y$ for all $k,l\in\mathbb{N}$,
		\item[(ii)] $L[E]$ is a bounded subset of $Y$, where $L(x)=Ax$ for all $x\in E$.
	\end{enumerate}
\end{lem}

\begin{proof}
Assume that $A \in (X, Y)$. Then $L[E]$ is bounded in $Y$, that is, $\|L(x)\|_Y \leq M \|x\|_X$ for all $x\in E$ and some positive constant $M$. For each fixed $k,l\in\mathbb{N}$, let's consider the $kl^{th}$ column of $A$, which is denoted by $A^{kl}$. We need to show that $A^{kl} \in Y$. 
Now, for each $k,l,m,n\in\mathbb{N}$, we have:
\begin{eqnarray*}
	A^{kl} = \sum_{m,n=1}^{\infty} a_{mnkl}e^{(mn)}.
\end{eqnarray*}
Now, if we consider the $mn^{th}$ term in the series defined by $A^{kl}$, we will have $a_{mnkl}$. Then, by the definition of $\vartheta$-convergence, 

\begin{eqnarray*}
A^{kl} = \vartheta-\lim_{m,n\to\infty}A^{kl}[m,n],
\end{eqnarray*}
where 
\begin{eqnarray*}
A^{kl}[m,n] = \sum_{i,j=1}^{m,n} a_{ijkl}e^{(ij)}.
\end{eqnarray*}
Since $A \in (X,Y)$ and $A^{kl}[m,n] = L[x^{[m,n]}]$, where $x^{[m,n]}$ is a sequence in $X$. Therefore, $A^{kl}[m,n] \in Y$ for all natural numbers $m$ and $n$, and by definition of $FDK$ space, the sequence $A^{kl}$ is convergent with respect to $\vartheta$ and hence belongs to $Y$. Thus, condition (i) is satisfied. Moreover, $A^{kl}[m,n] = L[x^{[m,n]}]$, and $L[E]$ is bounded in $Y$, it follows that $A^{kl}$ is bounded for all natural numbers $k$ and $l.$ Hence, condition (ii) suffices.

Conversely, assume that the conditions stated in (i) and (ii) hold true. Then, we have the following facts:
\begin{enumerate}
\item[(i)]  Since each column of $A$ belongs to $Y$, by definition of $(X,Y)$, we have $A^{kl} \in Y$ for all $k,l$.
\item[(ii)] Since $L[E]\subset Y$ is bounded, $\exists M>0$ satisfying $\|L(x)\|_Y \leq M \|x\|_X$ for all $x \in E$.
\end{enumerate}
Now, let $x \in X$ be arbitrary. By definition of $BDK$ space, $x = \vartheta-\lim_{m,n\to\infty}x^{[m,n]}$, where $x^{[m,n]} \in E$ for all natural numbers $m$ and $n$. Then, we have
\begin{eqnarray*}
	L(x) &= \vartheta-\lim_{m,n\to\infty}L[x^{[m,n]}] \\
	&= \vartheta-\lim_{m,n\to\infty}A^{[m,n]}.
\end{eqnarray*}
Since each $A^{kl}$ is bounded, it shows that $L(x)$ is also bounded, that is, $A \in (X,Y)$. This ends the proof.
\end{proof}

\begin{thm}\label{theorem4.6}
$A=(a_{mnkl})\in(\mathcal{H}_{\vartheta},\mathcal{H}_{\vartheta})$ iff 
\begin{eqnarray}
\label{thm4.6i}\vartheta-\lim_{m,n\to\infty}a_{mnkl}=0 \forall k,l\in\mathbb{N},\\
\label{thm4.6ii}\|A\|_{(\mathcal{H}_{\vartheta},\mathcal{H}_{\vartheta})}=\sup_{k,l\in\mathbb{N}}\left(\frac{1}{kl}\sum_{m,n=1}^{\infty}\left|\sum_{i,j=1}^{k,l}\Delta_{mn}^{11}a_{mnij}\right| \right),
\end{eqnarray}
where $\Delta_{mn}^{11}a_{mnij}=a_{mnij}-a_{m+1,nij}-a_{m,n+1,ij}+a_{m+1,n+1,ij}$
\end{thm}

\begin{proof}
Since $\mathcal{H}_{\vartheta}$ is a $BDK-$space with $AK(\vartheta)$ by Theorem \ref{thm1} and Theorem \ref{thm3}, and $E=\left\{\frac{1}{kl}e^{kl}:k,l\in\mathbb{N}\right\}$ is a determining set of $\mathcal{H}_{\vartheta}$ by Theorem \ref{thm4.4}, we apply Lemma \ref{Lemma 4.5} with $X=Y=\mathcal{H}_{\vartheta}$ which is a $BDK-$space with $AK(\vartheta)$ with its norm defined by (\ref{Hahn norm}). Let $k,l\in\mathbb{N}$ and $d^{(kl)}=\frac{1}{kl}e^{kl}$. Then 

\begin{eqnarray}
\label{eq4.6}A_{mn}d^{(kl)}=\sum_{i,j=1}^{\infty}a_{mnij}d_{ij}^{(kl)}=\frac{1}{kl}\sum_{i,j=1}^{k,l}a_{mnij}.
\end{eqnarray}
Therefore, we have two results as follows:
\begin{eqnarray}
\label{eq4.61} &&\|A_{mn}d^{(kl)}\|_{\mathcal{H}_{\vartheta}}<\infty \textit{ for all }d^{(kl)}\in E,\\
\label{eq4.62} &&A_{mn}d^{(kl)}\in\mathcal{C}_{\vartheta 0}\textit{ for all }d^{(kl)}\in E.
\end{eqnarray}
Then, we obtain by (\ref{eq4.61})

\begin{eqnarray*}
\|A_{mn}d^{(kl)}\|_{\mathcal{H}_{\vartheta}}&=&\sum_{m,n=1}^{\infty}mn\left|\Delta_{mn}^{11} A_{mn}d^{(kl)} \right|\\
&=&\sum_{m,n=1}^{\infty}mn\left|A_{mn}d^{(kl)}-A_{m+1,n}d^{(kl)}-A_{m,n+1}d^{(kl)}+A_{m+1,n+1}d^{(kl)} \right|\\
&=&\sum_{m,n=1}^{\infty}mn\left|\frac{1}{kl}\sum_{i,j=1}^{k,l}a_{mnij}-\frac{1}{kl}\sum_{i,j=1}^{k,l}a_{m+1,nij}-\frac{1}{kl}\sum_{i,j=1}^{k,l}a_{m,n+1,ij}+\frac{1}{kl}\sum_{i,j=1}^{k,l}a_{m+1,n+1,ij} \right|\\
&=&\frac{1}{kl}\sum_{m,n=1}^{\infty}mn\left|\sum_{i,j=1}^{k,l}\Delta_{mn}^{11}a_{mnij}\right|
\end{eqnarray*}
for all $k,l\in\mathbb{N}$. Thus it is clearly seen that the conditions in (\ref{eq4.61}) and (\ref{thm4.6ii}) are equivalent.

Moreover, the condition (\ref{eq4.62}) shows clearly that the condition in (\ref{thm4.6i}) also holds. This ends the proof.
\end{proof}

The given result is obtained by considering the determining set $E$ of the double sequence space $\mathcal{H}_{\vartheta}$ and by utilizing Lemma \ref{Lemma 4.5}. We employ the norm $\|A_{mn}d^{(kl)}\|_Y$ of the space $Y$, where $Y=\{\mathcal{BV}, \mathcal{BV}_{\vartheta 0}, \mathcal{CS}_{\vartheta},\mathcal{CS}_{\vartheta 0},\mathcal{BS}\}$, respectively.
\begin{cor}
Let $\vartheta=\{p,bp,r\}$. Then, each of the following statements is valid for 4D infinite matrix $A=(a_{mnkl})$.
\begin{enumerate}
\item[(i)] $A\in(\mathcal{H}_{\vartheta},\mathcal{BV})$ iff  
\begin{eqnarray}
\label{cor1.1}\sup_{k,l\in\mathbb{N}}\left(\frac{1}{kl}\sum_{m,n=0}^{\infty}\left|\sum_{i,j=1}^{k,l}\Delta_{mn}^{11}a_{mnij}\right|\right)<\infty.
\end{eqnarray}
\item[(ii)] $A\in(\mathcal{H}_{\vartheta},\mathcal{BV}_{\vartheta 0})$ iff  (\ref{thm4.6i}) and (\ref{cor1.1}) hold.
\item[(iii)] $A\in(\mathcal{H}_{\vartheta},\mathcal{CS}_{\vartheta})$ iff 
\begin{eqnarray}
\label{cor1.2}&& \sup_{s,t,k,l\in\mathbb{N}}\left(\frac{1}{kl}\sum_{m,n=0}^{s,t}\left|\sum_{i,j=1}^{k,l}a_{mnij}\right|\right)<\infty,\\
\label{cor1.3}&&\vartheta-\lim_{s,t\to\infty}\left(\frac{1}{kl}\sum_{m,n=0}^{s,t}\sum_{i,j=1}^{k,l}a_{mnij}\right)\textit{ exists for all }k,l\in\mathbb{N}.
\end{eqnarray}
\item[(iv)] $A\in(\mathcal{H}_{\vartheta},\mathcal{CS}_{\vartheta 0})$ iff  (\ref{cor1.2}) holds and 
\begin{eqnarray}
\label{cor1.4}\vartheta-\lim_{s,t\to\infty}\left(\frac{1}{kl}\sum_{m,n=0}^{s,t}\sum_{i,j=1}^{k,l}a_{mnij}\right)=0, \textit{ for all }k,l\in\mathbb{N}.
\end{eqnarray}
\item[(v)] $A\in(\mathcal{H}_{\vartheta},\mathcal{BS})$ iff  (\ref{cor1.2}) holds.
\end{enumerate}
\end{cor}

\begin{thm}\label{thm 4.80}
Let's establish the relationship between the terms of the  4D matrices $A=(a_{mnkl})$ and $F=(f_{mnkl})$ as follows:
\begin{eqnarray}\label{eq4-4}
	f_{mnkl}=\sum_{i,j=m-1,n-1}^{m,n}(-1)^{m+n-(i+j)}mna_{ijkl}=mn\Delta_{11}^{mn}a_{mnkl}.
\end{eqnarray}
Then, the 4D matrix $A=(a_{mnkl})\in(\lambda,\mathcal{H}_{\vartheta})$ iff  $F\in(\lambda,\mathcal{L}_u),$ where $\lambda$ is any given double sequence space.
\end{thm}

\begin{proof}
First, we assume that $A\in(\lambda,\mathcal{H}_{\vartheta})$. Therefore, for every $x\in\lambda$, $Ax$ exists, and contains in $\mathcal{H}_{\vartheta}$. Then, considering the relation stated in (\ref{eq4-4}), we deduce that
\begin{eqnarray*}
\sum_{k,l=0}^{s,t}f_{mnkl}x_{kl}=\sum_{k,l=0}^{s,t}\sum_{i,j=m-1,n-1}^{m,n}(-1)^{m+n-(i+j)}mna_{ijkl}x_{kl}=\sum_{i,j=m-1,n-1}^{m,n}(-1)^{m+n-(i+j)}mn\sum_{k,l=0}^{s,t}a_{ijkl}x_{kl}
\end{eqnarray*}
for all $m,n,s,t\in\mathbb{N}$. By letting $s,t\to\infty$, we deduce that
	
\begin{eqnarray*}
(Fx)_{mn}=\sum_{k,l}f_{mnkl}x_{kl}&=&\sum_{i,j=m-1,n-1}^{m,n}(-1)^{m+n-(i+j)}mn\sum_{k,l}a_{ijkl}x_{kl}\\
&=&mn\sum_{i,j=m-1,n-1}^{m,n}(-1)^{m+n-(i+j)}(Ax)_{ij}\\
&=&mn\Delta_{11}^{mn}(Ax)_{mn}
\end{eqnarray*}
for all natural numbers $m$ and $n$. Consequently, it can be observed that $Ax\in\mathcal{H}_{\vartheta}$ whenever $x\in\lambda$ iff  $Fx\in\mathcal{L}_u$ whenever $x\in\lambda$. Therefore, we conclude that $A\in(\lambda,\mathcal{H}_{\vartheta})$ iff  $F\in(\lambda,\mathcal{L}_{u})$.
	
\end{proof}

Let $\vartheta=\{bp,bp0,r,r0\}$. Since the proof can be done by taking the determining set $E$ of the double sequence space $\mathcal{H}_{\vartheta}$, the Lemma \ref{Lemma 4.5} and Theorem \ref{thm 4.80} into consideration, and by using the norm $\|A_{mn}d^{(kl)}\|_{\mathcal{H}_{\vartheta}}$ and the $BDK-$space $X$ with $AK(\vartheta)$, where $X=\{\mathcal{L}_u, \mathcal{C}_{\vartheta 0}, \mathcal{C}_{\vartheta},\mathcal{M}_{u}\}$, respectively, we give the following Theorem without its proof.
\begin{thm}
Assume $\vartheta=\{p,bp,r\}$. Then, each of the following assertions is valid for 4D infinite matrix $A=(a_{mnkl})$:
\begin{enumerate}
\item[(i)] $A\in(\mathcal{L}_u,\mathcal{H}_{\vartheta})$ iff  (\ref{thm4.6i}) holds and
\begin{equation}
	\label{cor2.1}\sup_{k,l\in\mathbb{N}}\sum_{m,n=1}^{\infty}mn\left|\Delta_{11}^{mn}a_{mnkl}\right|<\infty.
\end{equation}
\item[(ii)] $A\in(\mathcal{C}_{\vartheta 0},\mathcal{H}_{\vartheta})$ iff  (\ref{thm4.6i}) holds and
\begin{equation}
	\label{cor2.2}\sup_{K\subset\mathbb{N}~finite}\sum_{m,n=1}^{\infty}mn\left|\sum_{k,l\in K}\Delta_{11}^{mn}a_{mnkl}\right|<\infty.
\end{equation}
\item[(iii)] $A\in(\mathcal{C}_{\vartheta},\mathcal{H}_{\vartheta})$ iff  (\ref{thm4.6i}) and (\ref{cor2.2}) hold, and
\begin{equation}
	\label{cor2.3}\sup_{k,l\in\mathbb{N}}\sum_{m,n=1}^{\infty}mn\left|\sum_{i,j=1}^{k,l}\Delta_{11}^{mn}a_{mnij}\right|<\infty.
\end{equation}
\item[(iv)] $A\in(\mathcal{M}_{u},\mathcal{H}_{\vartheta})$ iff  
\begin{eqnarray}
\label{cor2.4}&&\sup_{K\subset\mathbb{N}~finite}\sum_{k,l=1}^{\infty}\left|\sum_{m,n\in\mathbb{N}}mn\Delta_{11}^{mn}a_{mnkl}\right|<\infty,\\
\label{cor2.5}&&\vartheta-\lim_{m,n\to\infty}\sum_{k,l=1}^{\infty}|a_{mnkl}|=0.
\end{eqnarray}

\end{enumerate}

\end{thm}

\section{Conclusion}
Most recently, Demiriz and Duyar \cite{Demiriz}, \c{C}apan and Ba\c{s}ar \cite{CHBF}, and Tu\v{g}, Rako\`{c}evi\`{c} and Malkowsky \cite{OVM3} introduced and studied the new spaces $\mathcal{M}_{u}(\Delta)$ and $\mathcal{BV}_q$, and $\mathcal{BV}_{\vartheta 0}$, where $\vartheta=\{p,bp,r\}$, respectively. Within this work, we introduced Hahn double sequence space $\mathcal{H}_{\vartheta}$, where $\vartheta=\{p,bp,r\}$ and we studied various topological and algebraic aspects, and we proved several containment relations. Moreover, we computed the algebraic dual and $\alpha-,\beta(bp)-,\gamma-$ duals of the Hahn double sequence space $\mathcal{H}_{\vartheta}$. Finally, we characterized the classes $(\mathcal{H}_{\vartheta},\lambda)$, where $\lambda=\{\mathcal{H}_{\vartheta},\mathcal{BV}, \mathcal{BV}_{\vartheta 0}, \mathcal{CS}_{\vartheta},\mathcal{CS}_{\vartheta 0},\mathcal{BS}\}$ and $(\mu,\mathcal{H}_{\vartheta})$, where $\mu=\{\mathcal{L}_u, \mathcal{C}_{\vartheta 0}, \mathcal{C}_{\vartheta},\mathcal{M}_{u}\}$.

As an extension of this research, building upon the foundational insights presented in the monographs \cite{Malkowsky and Rakocevic,Malafosse Malkowsky and Rakocevic,Banas Goebel}, there remain re-searchable problems concerning the determination of the Hausdorff measure of non-compactness and certain characteristics of matrix operators within the context of the Hahn double sequence space $\mathcal{H}_{\vartheta}$. Additionally, the computation of 4D matrix domains within the Hahn double sequence space $\mathcal{H}_{\vartheta}$, and characterization of some novel matrix classes from and into $\mathcal{H}_{\vartheta}$ remain open for further exploration.





\end{document}